\documentclass[10pt,showpacs,amsmath,amssymb,wasysym]{article}
\usepackage{units}
\usepackage{bbold,euler,newcent}
\usepackage{graphicx}
\usepackage{dcolumn}
\usepackage{bm,bbm,mathtools,esvect}
\usepackage{slashed,fullpage}

\title{ {\sc A primer on the differential geometry of quaternionic curves}\\ \vspace{5mm}}
\author{ {\tt SERGIO GIARDINO}\footnote{\tt sergio.giardino@ufrgs.br}\\ \vspace{0.5cm}\\
 Departamento de Matem\'atica Pura e Aplicada \\ Universidade Federal do Rio Grande do Sul (UFRGS)\\
Avenida Bento Gon\c calves 9500, Caixa Postal 15080, 91501-970\\  Porto Alegre, RS, Brazil}

\begin{document}
\date{} 

\newtheorem{theorem}{Theorem}[section]
\newtheorem{remark}{Remark}[section]
\newtheorem{lemma}{Lemma}[section]
\newtheorem{proposition}{Proposition}[section]
\newtheorem{corollary}{Corollary}[section]
\newtheorem{definition}{Definition}[section]

\maketitle

\begin{abstract}
\noindent 
This paper describes the foundations of a differential geometry of a quaternionic curves. The Frenet-Serret equations and the evolutes and evolvents of a particular quaternionic curve are accordingly determined. This new formulation takes benefit of the quaternionic structure and the results are much simpler than the present formulations of quaternionic curves.
\end{abstract}

\tableofcontents

\section{INTRODUCTION}

The study of quaternionic curves begins in \cite{Bharathi:1987kbn} where Frenet-Serret equations were determined. Several contributions have been made since then, and even the Frenet-Serret equations were further studied \cite{Gunes:1994mmt,Aksoyak:2019}, 
although different quaternic curves where this approach has been applied were considered as well	 \cite{Coken:2004otq,Hacisalihoglu:2011gok,Gungor:2011scq,Ilarslan:2013bma,EJPAM1979,Bektas:2016qoc,Kisi:2017aap,Coken:2017nqc,Karadag:2019scq,Kizilay:2020eqc,Kahraman:2020dqn}. These articles have the common feature that the Frenet-Serret
equations are strongly influenced by the geometric description, where a local basis is determined in each point of an arbitrary smooth and regular curve. 

In the present article, on the other hand, we obtain Frenet-Serret-like equations where the global algebraic structure of quaternions plays the most important role and effectively determines the equations. Using quaternions parametrized in polar coordinates, we develop two descriptions of quaternionic curves that can be understood as a projection of the quaternionic curve over a two dimensional plane.

To the best of our knowledge, this approach has never been reported, although two dimensional curves within a real space have a simple parametrization with complex numbers \cite{Hilario:2003dcg}, and here we closely follow this approach. The example of evolutes and evolvents presented here shows that our proposal is much simpler than the current approach \cite{Gungor:2013qie,Onder:2020gqi,Altun:2017scy,Senyurt:2017scy,Hou:2018gia}, based in \cite{Bharathi:1987kbn}.  Before developing the differential geometric approach, in the next section we present the most important features of quaternions that will be used in this article.

\section{QUATERNIONIC BACKGROUND}
Quaternions ($\mathbbm H$) are hyper-complex numbers that, if $\,q\in\mathbbm H$, then
\begin{equation}\label{hn01}
q=x_0+x_1i+x_2 j+x_3k,\qquad\mbox{where}\qquad x_\mu\in\mathbbm R\qquad\mbox{for}\qquad\mu\in\{0,\,1,\,2,\,3\},
\end{equation} 
and $\,i,\,j\,$ and $\,k\,$ are anti-commuting imaginary units that satisfy
\begin{equation}\label{hn02}
i^2\,=\,j^2\,=\,k^2\,=-1\qquad\mbox{and}\qquad ijk=-1.
\end{equation} 
In the same fashion as complex numbers, the quaternionic conjugate and the quaternionic norm are
\begin{equation}\label{hn03}
\overline q\,=\,x_0-x_1i-x_2 j-x_3k,\qquad\qquad |q|^2\,=\,q\overline q\,=\,x_0^2+x_1^2+x_2^2+x_3^2.
\end{equation} 
According to Adolf Hurwitz \cite{Hurwitz:1898hqt}, quaternions comprise one of the four division algebras, and the other are the reals ($\mathbbm R$), the complexes ($\mathbbm C$) and the octonions ($\mathbbm O$). There are many sources that introduce and describe the quaternions \cite{Ward:1997qcn,Garling:2011zz,Rocha:2013qtt}, and in this section we focus on specific notations for  representing quaternions that will be useful in this article. The notation (\ref{hn01}) can be named the extended or Cartesian notation, which additionally may be written as
\begin{equation}\label{hn04}
q=x_0+\omega |\bm x|\qquad \mbox{where}\qquad |\bm x|^2=x_1^2+x_2^2+x_3^2\qquad \mbox{and}\qquad\omega=\frac{\bm x}{|\bm x|}.
\end{equation}
This notation admits a complex unit, so that $\omega^2=-1$, and $x_0$ and $\bm x$ are respectively named the scalar (or temporal) and the vector (or spatial) components of the quaternion. This notation has the convenient property that the scalar and vector  components are commutative, and the liability that every quaternionic number has its own imaginary unit $\omega$ that neither commutes nor anti-commutes with any other quaternionic unit.
We also observe that the imaginary component of (\ref{hn04}) is always positive. Using this feature, we obtain the polar notation of the Cartesian quaternion
\begin{equation}\label{hn05}
q\,=\,\rho\big(cos\theta+\sin\theta\,\omega\big)\qquad\mbox{where}\qquad\rho=|q|\qquad\mbox{and}\qquad\theta\in[0,\,\pi].
\end{equation}
In addition to the Cartesian notation (\ref{hn01}), in symplectic notation the quaternionic number reads
\begin{equation}\label{hn06}
q=z_0+z_1 j\qquad\mbox{where}\qquad z_0=x_0+x_1 i\qquad\mbox{and}\qquad z_1=x_2 +x_3 i.
\end{equation}
The symplectic notation is not unique, and we can replace (\ref{hn06}) with 
\begin{equation}\label{hn07}
q=z_0+\overline\zeta\, k \qquad\mbox{where}\qquad \zeta=x_3 +x_2 i,
\end{equation}
and four other possibilities replacing $i$ with $j$ and $k$ in (\ref{hn06}). A symplectic polar form is 
\begin{equation}\label{hn08}
q\,=\,\rho\Big(\cos\vartheta e^{i\phi}+\sin\vartheta e^{i\psi}j\Big)\qquad\mbox{where}\qquad \vartheta\in\left[0,\frac{\pi}{2}\right]\qquad\mbox{and}
\qquad \phi,\,\psi\in[0,\,2\pi].
\end{equation}\label{hn09}
From \cite{Harvey:1990sca}, we adopt the scalar product for quaternions
\begin{equation}\label{hn10}
\big\langle p,\,q\big\rangle\,=\,\mathfrak{Re}\big[\,p\overline q\,\big],
\end{equation}
which permit us to define that if $p$ and $q$ are orthogonal, then
\begin{equation}\label{hn11}
\big\langle p,\,q\big\rangle\,=\,0
\end{equation}
and if
\begin{equation}\label{hn12}
\big\langle p,\,q\big\rangle\,=\,\,p\overline q,
\end{equation}
then $p$ and $q$ are parallel. Properties (\ref{hn10}-\ref{hn12}) are also valid for complex numbers. 
Using this definition, we can obtain the orthogonality relations. In terms of the Cartesian notation,
\begin{equation}\label{hn13}
\langle q,\,e_i q\rangle = 0\qquad\mbox{where}\qquad e_i=\{i,\,j,\,k\}.
\end{equation}
and thus we have that quaternions can be described as a four dimensional real vector space. In the polar notations (\ref{hn05}) and (\ref{hn06}), we respectively obtain, 
\begin{equation}\label{hn14}
\langle q,\,\omega q\rangle =0,\qquad\mbox{and}\qquad\langle q,\, q j\rangle =0,
\end{equation}
and thus the quaternions are described as two dimensional real vector spaces. A strange thing, if we consider that the space had originally four dimensions from (\ref{hn13}). In order to clarify this point, let us remember that 
 the restriction of the range of $\theta$ in (\ref{hn05}) and $\vartheta$ in (\ref{hn08}) indicates that a negative signal has to be absorbed by a complex structure in order to keep the polar angles within the correct range.
Let us consider the simplest case, the negative of a Cartesian polar quaternion
\begin{equation}\label{hn15}
q=(-1)(\cos\theta+\omega\sin\theta)
\end{equation}
If $q$ were a complex number, this multiplication means a rotation of the polar angle, so that $\theta\to\theta+\pi$, something prohibited in the quaternionic case. Thus, the negative signal must be absorbed by the imaginary unit $\omega$, and we the correct rotation will be $\theta\to\pi-\theta$. Furthermore, this rotation changes the orientation of the polar angle, and hence the quaternionic polar angle is not orientable in this case. A wishful result in dimensions higher than two, where no orientable polar angle is expected. Using the above principle, we obtain
\begin{equation}\label{hn16}
\left.
\begin{array}{l}
\quad	\;\, q\,=\,\cos\theta+\omega\sin\theta\\ \\
\;\;\omega q\,=\,\cos\left(\theta+\frac{\pi}{2}\right)+\omega\sin\left(\theta+\frac{\pi}{2}\right)\\ \\
\omega^2 q\,=\,\cos\left(\pi-\theta\right)+\widetilde\omega\sin\left(\pi-\theta\right)\\ \\
\omega^3 q\,=\,\cos\left(\frac{\pi}{2}-\theta\right)+\widetilde\omega\sin\left(\frac{\pi}{2}-\theta\right)\\
\end{array}
\right\}
\qquad\theta\in \left[0,\frac{\pi}{2}\right]\qquad\mbox{and}\qquad \widetilde\omega=-\omega
\end{equation}
The case for $\,\theta\in\left[\frac{\pi}{2},\,\pi\right]\,$ can also be easily obtained, and hence the conditions of (\ref{hn05}) are satisfied in every case. We can also  understand that we have two orthogonal pairs $(q,\,\omega q)$ and $(\omega^2 q,\,\omega^3 q)$, 
thus counting the four components for the space, which are separated because their imaginary unit is different. Using this ideas, in general we have
\begin{equation}\label{hn17}
q(\theta_1+\theta_2)\,=\,\left\{
\begin{array}{ll}
\cos\theta_0+\omega\sin\theta_0 & n\quad\mbox{even}\qquad \theta_0\in\left[0,\,\pi\right]\\ \\
\cos\left(\pi-\theta_0\right)+\widetilde\omega\sin\left(\pi-\theta_0\right)& n\quad\mbox{odd,}\qquad \theta_0\in\left[0,\,\frac{\pi}{2}\right]\\ \\
\cos\left(\frac{\pi}{2}-\theta_0\right)+\widetilde\omega\sin\left(\frac{\pi}{2}-\theta_0\right)
& n\quad\mbox{odd,}\qquad \theta_0\in\left[\frac{\pi}{2},\,\pi\right]\\ \\
\quad\;\,\theta_1+\theta_2\,=\,n\pi+\theta_0\qquad \mbox{and}\qquad \widetilde\omega=-\omega&
\end{array}
\right.
\end{equation}
The above analysis can be repeated to the symplectic case (\ref{hn08}), so that
\begin{equation}\label{hn18}
\left.
\begin{array}{l}
\;\;\; q\, =\,\cos\vartheta e^{i\phi}+\sin\vartheta e^{i\psi}j\\ \\
\;\, q j\,=\,\cos\left(\frac{\pi}{2}-\vartheta\right)e^{i(\psi-\pi)}+\sin\left(\frac{\pi}{2}-\vartheta\right)e^{i\phi}j\\ \\
q j^2\,=\,\cos\vartheta e^{i(\phi-\pi)}+\sin\vartheta e^{i(\psi-\pi)}j\\ \\
q j^3\,=\,\cos\vartheta\left(\frac{\pi}{2}-\vartheta\right)e^{i\psi}+\sin\left(\frac{\pi}{2}-\vartheta\right)e^{i(\phi-\pi)}j\\
\end{array}
\right\}
\qquad\vartheta\in \left[0,\frac{\pi}{2}\right]
\end{equation}
and also
\begin{equation}\label{hn19}
q(\vartheta_1+\vartheta_2)\,=\,
\left\{
\begin{array}{ll}
\cos\vartheta_0 e^{i\phi}+\sin\vartheta_0 e^{i\psi}j & \quad n=0\;\mbox{mod}\;4\\ \\
\cos\left(\frac{\pi}{2}-\vartheta_0\right)e^{i(\psi-\pi)}+\sin\left(\frac{\pi}{2}-\vartheta_0\right)e^{i\phi}j & \quad n=1\;\mbox{mod}\;4\\ \\
\cos\vartheta_0 e^{i(\psi-\pi)}+\sin\vartheta_0 e^{i(\phi-\pi)}j &\quad n=2\;\mbox{mod}\;4\\ \\
\cos\vartheta\left(\frac{\pi}{2}-\vartheta_0\right)e^{i\psi}+\sin\left(\frac{\pi}{2}-\vartheta_0\right)e^{i(\phi-\pi)}j&\quad n=3\;\mbox{mod}\;4
\end{array}
\right.
\end{equation}
where
\[
\vartheta_1+\vartheta_2=\vartheta_0 + n\frac{\pi}{2},\qquad\vartheta_0\in\left[0,\,\frac{\pi}{2}\right]
\qquad\mbox{and}\qquad n\in\mathbbm N.
\]
Now we have every necessary concept to develop a differential geometry of quaternionic curves: the polar notation, the scalar product and the orthogonality condition. We develop this in the next section.
\section{QUATERNIONIC CURVES\label{QC}}
In this section we present original results concerning the description of quaternionic curves. We divide this section into two subsections, in the first one we develop the formalism of Frenet-Serret using the polar notation of the Cartesian quaternion (\ref{hn05}), and in the second subsection we repeat the procedure using the symplectic parametrization (\ref{hn06}). Although the results are similar, it is important to consider them separately because each quaternionic curve may have a more suitable description, or even an unique description, using a particular parametrization.

\subsection{polar quaternion curves}
We define define a quaternic curve using a standard way
\begin{definition}\label{qc01} A curve in $\mathbbm H$ is the smooth application
$q(t): I\to\mathbbm H\,$ for $\,I\subset\mathbbm R\,$  as
\[
q=x_0+x_1i+x_2 j+x_3k,
\]
where $\,x_\mu=x_\mu(t)$ are real functions $\mathcal C^\infty$ for $\,t\in I\,\mbox{ and }\,\mu\in\{0,\,1,\,2,\,3\}.\,$ 
\end{definition}
The tangent application $\,T\,$ of a quaternionic curve $\,q\,$ is
\begin{equation}\label{qc02}
T=\frac{dq}{dt}=q'.
\end{equation}
The length of a quaternionic curve is also standardly defined as the real function 
\begin{equation}\label{qc03}
\mathcal L(t)=\int_{t_0}^t\big|T(t')\big|dt'.
\end{equation}
It is also standard that a curve can be parametrized by their length, and that a curve is parametrized by its length if and only if
$\,T(t)=1.\,$ As a last concept, we define:
\begin{definition} 
A quaternionic curve $q:\mathbbm I\to\mathbbm H$ is regular in $\,t_0\in I\,$ if $\,T(t_0)\neq 0,\,$ and is regular in $\,I\,$ if it is regular to $\,\forall\, t\in I.\,$
\end{definition}
Now, we have every necessary element to the description of quaternionic curves. From the definition of scalar product (\ref{hn10}) and from (\ref{hn13}),  $\,iq,\,jq\mbox{ and } kq\,$ are orthogonal to $q$. Thus, there are three applications $N_i$ that are normal to the tangent application $T$, so that
\begin{equation}\label{qc04}
N_i\,=\,e_i T\qquad\mbox{and evidently}\qquad \big\langle N_i,\,T\big\rangle =0.
\end{equation}
Supposing that $\,T\,$ is unitary, we obtain
\begin{equation}\label{qc05}
\left(T \overline{T}\right)'\,=0\qquad\Rightarrow\qquad T'\overline T+T\overline T'=\big\langle T,\,T'\big\rangle = 0
\end{equation}
and thus $T$ and $T'$ are orthogonal. Hence $\,T'\,$ belongs to the three dimensional space generated by (\ref{qc04}), and can be written as
\begin{equation}\label{qc06}
T'\,=\,\sum_i \kappa_i N_i
\end{equation}
where $\,\kappa_i=\kappa_i(t)\,$ are real functions. Equation (\ref{qc06}) can be written as
\begin{equation}\label{qc07}
q''\,=\,\kappa q',\qquad\mbox{where}\qquad \kappa=\kappa_1 i+\kappa_2 j+\kappa_3 k.
\end{equation}
The pure imaginary quaternion $\kappa$ may be called the quaternionic curvature.
The quaternionic equation (\ref{qc07}) can be turned in four real equations, namely
\begin{equation}\label{qc08}
\left(
\begin{array}{c}
x_0 \\
x_1 \\
x_2 \\
x_3 \\
\end{array}
\right)''\,=\,
\left(
\begin{array}{cccc}
0 & -\kappa_1 & -\kappa_2 & -\kappa_3 \\
\kappa_1 & 0 & -\kappa_3 & \kappa_2 \\
\kappa_2 & \kappa_3 & 0 & -\kappa_1 \\
\kappa_3 & -\kappa_2 & \kappa_2 & 0 \\
\end{array}
\right)
\left(
\begin{array}{c}
x_0 \\
x_1 \\
x_2 \\
x_3 \\
\end{array}
\right)'.
\end{equation}
Equations (\ref{qc08}-\ref{qc08}) are equivalent to Frenet-Serret equations, and are different from the current formulation 
\cite{Bharathi:1987kbn,Gunes:1994mmt,Aksoyak:2019}. However, differently to this current formulation, where the basis changes at every point of the space, in this case the basis is determined by the quaternic directions. We point out that the matrix of equations (\ref{qc08}) agrees with the matrix formulation of quaternions \cite{Ward:1997qcn}. In the case of a non-unitary $|q'|,$ we write
\begin{equation}\label{qc09}
T=\frac{q'}{|q'|}.
\end{equation}
$T$ is obviously parallel to $q',$ and consequently their orthogonal space are common. Consequently
\begin{proposition}\label{qcP01} For a regular $q\in\mathbb H$ holds that
\[
\frac{d}{dt}\left(\frac{q'}{|q'|}\right)\,=\,\kappa q'.
\]
\end{proposition}
$\hfill\triangle$

This proposition permits us to state:
\begin{proposition}\label{qcP02} For a regular $q\in\mathbbm H,$ it holds that
\[
\kappa_i=\frac{1}{|q'|^2}\big\langle q'',\,N_i\big\rangle\qquad where \qquad N_i=e_i\frac{q'}{|q'|}.
\]
\rm{
{\bf Proof:} From Proposition \ref{qcP01} we have
\[
\frac{q''}{|q'|}+\left(\frac{1}{|q'|}\right)q'=|q'|\Big(\kappa_1 N_1+\kappa_2 N_2+\kappa_3 N_3\Big)
\]
Using $\,\langle N_i,\,T\rangle=0\,$ and $\,\langle N_i,\,N_j\rangle=\delta_{ij},\,$ we immediately obtain the result.
}
\end{proposition}
$\hfill\triangle$

The simplest solution for (\ref{qc07}) involving a constant quaternic curvature and may easily solved  using the polar parametrization 
(\ref{hn05}), so that
\begin{equation}\label{qc10}
q(t)=\frac{1}{|\kappa|}\Bigg[\cos\Big(|\kappa|t+\phi_0\Big)+\omega \sin\Big(|\kappa|t+\phi_0\Big)\Bigg],\qquad\mbox{where}\qquad \omega \,=\,\frac{\kappa}{|\kappa|},
\end{equation}
and $\phi_0$ is an arbitrary real phase. Using the above structure, many results of the differential geometry of plane curves can be proven.

\begin{theorem}[Fundamental theorem of quaternic curves]\label{qcT01} Let $|\kappa|:I\subset \mathbbm R\to\mathbbm R$ be a $\mathcal C^\infty$ function, $t_0\in I$, $P_0,\,V_0\in \mathbb H$ and $|V_0|=1$. Hence there is a unique curve $\,q(t):I\to\mathbbm H\,$ parametrized by the arc length whose curvature at each point is $|\kappa(t)|$, their initial value is $q(t_0)=P_0$ and their initial first derivative is $q'(t_0)=V_0$.

\rm{
{\bf Proof:} Equation (\ref{qc07}) can be solved using $\,\omega'=0\,$ and therefore
\[
q(t)=x(t)+\frac{\kappa(t)}{|\kappa(t)|}y(t)\qquad \mbox{where} \qquad
\left\{
\begin{array}{l}
x(t)=x_0\,+\int_{t_0}^{t''}\cos\left(\int_{t_0}^t|\kappa(t')|dt' +\phi_0\right)dt''\\
\\
y(t)=y_0\,+\,\int_{t_0}^{t''}\sin\left(\int_{t_0}^t|\kappa(t')|dt' +\phi_0\right)dt''
\end{array}\right. .
\]
The initial conditions are fixed as
\[
P_0\,=\,x_0+\cos\phi_0+\frac{\kappa(t_0)}{|\kappa(t_0)|}\Big(y_0+\sin\phi_0\Big)\qquad\mbox{and}\qquad 
V_0\,=\,-\sin\phi_0+\frac{\kappa(t_0)}{|\kappa(t_0)|}\cos\phi_0.
\]
If a certain function $p$ satisfies the conditions of the theorem, thus $\Delta=p-q$ satisfies (\ref{qc07}) for $\Delta(t_0)=\Delta'(t_0)=0$. The equation written in components implies that $|\Delta(t)|'=0$ for every $t$, and thus $p'=q'.$ As $\Delta(t_0)=0$, thus $\,p=q,\,\forall\,t\in I\,$  and the solution is unique.
}
\end{theorem}
$\hfill\triangle$

The results indicate that an unique quaternionic curve can uniquely be ascribed to an quaternionic curvature, and the correspondence to a complex parametrization of $\,\mathbbm R^2\,$ is exact, a wishful result that will be proved in the symplectic notation as follows.
\subsection{quaternion symplectic curves}
In this case, we use the quaternic notation (\ref{qc08}), and the quaternic Frenet-Serret equation (\ref{qc07}) turns into
\begin{equation}\label{qc12}
q''=q'c j\qquad\mbox{where}\qquad c: I\to\mathbb C
\end{equation}
is the complex curvature. In matrix form, (\ref{qc12}) is
\begin{equation}\label{qc14}
\left( \begin{array}{c} z_0\\z_1\end{array}\right)''=
\left( \begin{array}{cc} 0 & -\,\overline c\\c& 0\end{array}\right)
\left( \begin{array}{c} z_0\\z_1\end{array}\right)',
\end{equation}
and Proposition (\ref{qcP01}) becomes
\begin{equation}\label{qc15}
\frac{d}{dt}\left(\frac{q'}{|q'|}\right)\,=\, q'cj.
\end{equation}
The curvature $\,c\,$ is obtained in the same way as in Proposition (\ref{qcP02}), so that 
\begin{equation}\label{qc16}
\big\langle N\,\overline c,\,N\big\rangle\,=\,\frac{1}{|q'|^2}\big\langle q'',\,N_i\big\rangle\qquad \mbox{where}\qquad N=\frac{q'}{|q'|}j.
\end{equation}
However, $\,\big\langle N\,\overline c,\,N\big\rangle\,=\,\mbox{Re}[c],\,$ and the imaginary component of the curvature is obtained from
\begin{equation}\label{qc17}
\big\langle N\,\overline c\, i,\,N\big\rangle\,=\,\frac{1}{|q'|^2}\big\langle q''i,\,N\big\rangle 
\end{equation}
The simplest solution of (\ref{qc12}) is
\begin{equation}\label{qc18}
q(t)=\frac{1}{|c|}\Bigg[\cos\Big(|c|t+\phi_0\Big)e^{-i\frac{\theta}{2}}\,+\,\sin\Big(|c|t+\phi_0\Big)e^{i\frac{\theta}{2}}j\Bigg],\qquad\mbox{where}\qquad c=|c|e^{i\theta},
\end{equation}
On the other hand,  equation (\ref{qc12}) admits the general solution
\begin{equation}\label{qc19}
z'_0\,=\,e^{i\int_{t_0}^t|c(t')|dt'},\qquad z'_1\,=\,z'_0\, e^{i\phi_0}\qquad \mbox{and}\qquad
c\,=\,|c|e^{i\left(\phi_0+\frac{\pi}{2}\right)}.
\end{equation}
where $\phi_0$ is a real constant. Evidently,
\begin{equation}\label{qc20}
z_0(t)\,=\,P_0\,+\,\int_{t_0}^t z_0(t')dt',\qquad z_1(t)\,=\,Q_0\,+\,e^{i\phi_0}\int z_0(t') dt'
\end{equation}
where $P_0$ and $Q_0$ are complex constants and Theorem \ref{qcT01} also holds.

\section{EVOLUTES AND EVOLVENTS}
There are several former studies on this theme \cite{Gungor:2013qie,Onder:2020gqi,Altun:2017scy,Senyurt:2017scy,Hou:2018gia},
and their development is much more complicated than ours. In this section, we provide simple solutions for these solutions using a generalization to quaternionic curves that are analogous to the complex formulation for plane real curves.

\subsection{evolutes}

We define equation for the center of curvature as
\begin{equation}\label{qc21}
q_E(t)\,=\,q(t)+\omega\frac{1}{\,|\kappa|}q'(t),
\end{equation}
where the regular quaternic curve $\,q(t)\,$ is parametrized with their length and that satisfies (\ref{qc07}). Their first derivative is
\begin{equation}\label{qc22}
q'_E\,=\,\omega\left(\frac{1}{|\kappa|}\right)' q',
\end{equation}
and consequently $q'_E$ is normal to $q'$. Immediately, we have the proposition:
\begin{proposition} \label{qcP03} Given a regular curve $q$ of quaternic curvature $\kappa$ parametrized by the curve length, the elements $\kappa_{Ei}$ of  quaternionic curvature elements of their curvature center $q_E$ are given by
\[
\kappa_{Ei}\big|q'_E\big|=\kappa_i
\]
\end{proposition}
{\bf Proof:} The result immediately comes after calculating the element of $\kappa_E$  from Proposition (\ref{qcP01}) using (\ref{qc07}) and $|q'|=1$.

$\hfill\triangle$

\vspace{1mm}

The results (\ref{qc21}-\ref{qc22}) and Proposition \ref{qcP03} are  analogously obtained using the complex parametrization of plane curves, and thus we have that (\ref{qc21}) comprises the evolute of the quaternionic curve. In symplectic notation, the center of curvature and their first derivative are
\begin{equation}\label{qc23}
q_E(t)\,=\,q(t)+q'(t)\frac{c}{\,|c|^2}j\qquad\mbox{and}\qquad q'_E\,=\,q'\left(\frac{1}{|c|}\right)'\frac{c}{|c|}j
\end{equation}
and consequently $q_E'$ is normal to $q'$. Finally, we have
\begin{equation}\label{qc24}
\mathfrak{Re}[c_E]\big|q'_E\big|=\mathfrak{Re}[c]\qquad\mbox{and}\qquad\mathfrak{Im}[c_E]\big|q'_E\big|=\mathfrak{Im}[c]
\end{equation}
And the correspondence of the extended cases with the complex case is also exact.
\subsection{evolvents}

In the case of the evolvent (or involute) of a regular quaternionic curve $\,q,\,$ of unitary tangent application ($|q'(t)|=1$), we have
\begin{equation}\label{qc25}
q_I(t)=q(t)+\Big(\lambda_0-\mathcal L(t)\Big) q'(t)
\end{equation}
where $\,\lambda_0\,$ is a real constant and $\,\mathcal L(t)\,$ is the length function (\ref{qc03}). The first derivative
reads
\begin{equation}\label{qc26}
q'_I(t)\,=\,\Big(\lambda_0-\mathcal L(t)\Big)\kappa q'
\end{equation}
where the Proposition (\ref{qcP01}) were used. Again, we have that the tangent application of the evolvent is orthogonal to the tangent application of the curve. Using (\ref{qc07}) for $\,q_I(t)\,$, we obtain that the quaternionic curvature of the evolvent is
we obtain
\begin{equation}\label{qc27}
\kappa_I(t)=\frac{1}{\big|\lambda_0-\mathcal L(t)\big|}\omega,
\end{equation}
in sound agreement to the complex case. Finally, we obtain our last result:
\begin{proposition}
The quaternionic regular curve $\,q\,$ is the evolute of any of its evolvents.
\end{proposition}
\[
q(t)\,=\,\big(q_i\big)_E(t)
\]
The proof is obtained by immediate substitution of (\ref{qc25}) into (\ref{qc21})

In symplectic notation we can obtain every result obtained for envolvents description of (\ref{qc25}), we have that
\begin{equation}
q'_I(t)\,=\,\Big(\lambda_0-\mathcal L(t)\Big) q'\,c\,j,
\end{equation}
and the analogue of (\ref{qc27}) is
\begin{equation}
c_I=\frac{1}{\big|\lambda_0-\mathcal L(t)\big|}\frac{c}{|c|}
\end{equation}

This last result shows the complete agreement betwixt the description of quaternionic curves with the complex formulation of plane curves.
\section{CONCLUSION}
In this article we presented a novel description of quaternionic curves, where the most important role is played by the 
global algebraic structure of the quaternions. The equations are much simpler than the usual equations of the previous
cases, and we obtained general solutions using the polar system of coordinates for the Cartesian extended notation and also the symplectic representation of quaternions. The results generalize the two dimensional complex case, and we expect that this formulation will enable to shed brighter light over the previous results of quaternionic curves. Furthermore, we expect that additional results of plane curves can be transposed to four dimensions in order to clarify whether additional two dimensional cases could be generalized using the quaternionic formalism.
%
%
%
%

\bibliographystyle{unsrt} 
\bibliography{bib_gdh}

\end{document}